\documentclass[a4paper,abstracton]{scrartcl}

\usepackage{hyperref}

\usepackage[utf8]{inputenc}

\usepackage{amsmath}
\usepackage{amsthm}
\usepackage{amssymb}

\usepackage[usenames,dvipsnames]{color}
\usepackage{graphicx}
\usepackage{caption, subcaption}

\usepackage{tikz}
\usetikzlibrary{automata}
\usetikzlibrary{arrows,shapes}

\usepackage{bbding}

\usepackage{amssymb,amsmath,verbatim,tabularx,enumitem,colonequals,graphicx,psfrag}

\newcommand{\cU}{\mathcal{U}}
\newcommand{\X}{\mathcal{X}}
\newcommand{\R}{\mathbb{R}}
\newcommand{\cR}{\mathcal{R}}

\usepackage{authblk}
\usepackage[normalem]{ulem} 

\usepackage{framed}
\usepackage{rotating}
\usepackage{multirow}

\allowdisplaybreaks

\begin{document}

\title{An Experimental Comparison of Uncertainty Sets for Robust Shortest Path Problems}

\author{Trivikram Dokka}
\author{Marc Goerigk}
\affil{Department of Management Science, Lancaster University, Lancaster, United Kingdom}

\date{}

\maketitle

\begin{abstract}
Through the development of efficient algorithms, data structures and preprocessing techniques, 
real-world shortest path problems in street networks are now very fast to solve. But in reality, the exact travel times along each arc in the network may not be known. This lead to the development of robust shortest path problems, where all possible arc travel times are contained in a so-called uncertainty set of possible outcomes.

Research in robust shortest path problems typically assumes this set to be given, and provides complexity results as well as algorithms depending on its shape. However, what can actually be observed in real-world problems are only discrete raw data points. The shape of the uncertainty is already a modelling assumption. In this paper we test several of the most widely used assumptions on the uncertainty set using real-world traffic measurements provided by the City of Chicago. We calculate the resulting different robust solutions, and evaluate which uncertainty approach is actually reasonable for our data. This anchors theoretical research in a real-world application and allows us to point out which robust models should be the future focus of algorithmic development.
\end{abstract}


\section{Introduction}

The problem of finding shortest paths in real-world networks has seen considerable algorithmic improvements over the last decade \cite{bast2016route}. In the typical problem setup, one assumes that all data is given exactly. But also robust shortest path problems have been considered, where travel times are assumed to be given by a set of possible scenarios. In \cite{yu1998robust}, it was shown that the problem of finding a path that minimizes the worst-case length over two scenarios is already weakly NP-hard. For general surveys on results in robust discrete optimization, we refer to \cite{aissi2009min,kasperski2016robust}.

There are many possibilities how to model the scenario set that is used for the robust optimization process (see, e.g., \cite{montemanni2004exact,busing2012recoverable}), and it is not obvious which is ''the right'' one. Part of the current research ignores the problem by simply assuming that the uncertainty was ''given'' in some specific form, while this does not happen in reality.

In fact, the starting point for all uncertainty sets is raw data, given as a set of observations of travel times. This data is then processed to fit different assumptions on the shape and size of the uncertainty set, and preferences of the decision maker. So far, the discussion of these uncertainty sets has been led by theoretical properties, such as the computational tractability of the resulting robust model. We believe that this leads to a gap in the literature, where models are not sufficiently underpinned by actual real-world data to verify results. The purpose of this paper is to close this gap. We use real-world traffic observations by the City of Chicago to create a selection of the best-known and most-used uncertainty sets from the research literature. Using these uncertainty sets, we calculate different robust solutions and compare their performance. This allows us to determine which uncertainty sets are actually valuable for real-world robust shortest path problems. Our results give strong impulses for future research in the field by pointing out which problems are most worthy to solve more efficiently.

In Section~\ref{sec:unc} we briefly introduce all six uncertainty sets used in this study, and discuss the complexity of the resulting robust problems. The experimental setup and results are then presented in Section~\ref{sec:exp}, before concluding this paper in Section~\ref{sec:conclusion}.

\section{Uncertainty Sets for the Shortest Path Problem}
\label{sec:unc}

Let a directed graph $G=(V,A)$ with nodes $V$ and arcs $A$ be given. In the classic shortest path problem, each arc $e$ has some specific travel time $c_e \ge 0$. Given a start node $s$ and a target node $t$, the aim is to find a path minimizing the total travel time, i.e., to solve
\[ \min \{ \pmb{c}\pmb{x} : \pmb{x}\in\X \} \]
where $\X\subseteq\{0,1\}^n$ denotes the set of $s$-$t$-paths, and $n = |A|$. For our setting we assume instead that a set $\cR$ of travel time observations is given, $\cR = \{\pmb{c}^1,\ldots,\pmb{c}^N\}$ with $\pmb{c}^i\in\R^n$. This is the available raw data. In the well-known robust shortest path problem we assume that an uncertainty set $\cU$ is produced based on this raw data, and solve the robustified problem
\[ \min_{\pmb{x}\in\X} \max_{\pmb{c}\in\cU} \pmb{c}\pmb{x}, \]
that is, we search for a path that minimizes the worst-case costs over all scenarios. In the following sections we detail different possibilities from the current literature to generate $\cU$ based on $\cR$. Each set is equipped with a scaling parameter to control its size.

\subsection{Convex Hull}

In this approach, also known as discrete uncertainty (see \cite{yu1998robust,kasperski2016robust}), we set $\cU^{CH} = \cR$. The resulting robust problem can then be written as
\begin{align*}
\min \ & z \\
\text{s.t. } & z \ge \pmb{c}^i\pmb{x} & \forall i\in[N] \\
& \pmb{x}\in\X
\end{align*}
Note that this problem is equivalent to using $\cU^{CH} = conv (\{\pmb{c}^1,\ldots,\pmb{c}^N\})$. The problem is known to be NP-hard already for two scenarios. 

\noindent\textbf{Scaling:} 
Let $\hat{\pmb{c}}$ be the average of $\{\pmb{c}^1,\ldots,\pmb{c}^N\}$, i.e., $\hat{\pmb{c}} = \frac{1}{N}\sum_{i\in[N]} \pmb{c}^i$. For a given $\lambda\ge 0$, we substitute each point $\pmb{c}^i$ with $\hat{\pmb{c}} + \lambda(\pmb{c}^i-\hat{\pmb{c}})$, and take the convex hull of the scaled data points.

\subsection{Intervals}

We set $\cU^I$ as the smallest hypercube containing all data, i.e., $\cU^I = \prod_{i\in[n]} [\min_{j\in[N]} c^j_i,\max_{j\in[N]} c^j_i]$. For ease of notation, we write $\overline{c}_i := \max_{j\in[N]} c^j_i$ and $\underline{c}_i := \min_{j\in[N]} c^j_i$. The resulting robust problem is then
\begin{align*}
\min \ &\overline{\pmb{c}}\pmb{x} \\
\text{s.t. } & \pmb{x}\in\X
\end{align*}
which is a classic shortest path problem. Robust shortest path problems with interval uncertainty are therefore easy to solve, but frequently used, especially in the so-called min-max regret setting (see \cite{chassein2015new}).


\noindent\textbf{Scaling:} We use $\cU^I = \prod_{i\in[n]} \left[ 
\frac{\overline{c}_i + \underline{c}_i}{2} - \lambda \frac{\overline{c}_i - \underline{c}_i}{2},
\frac{\overline{c}_i + \underline{c}_i}{2} + \lambda \frac{\overline{c}_i - \underline{c}_i}{2}\right]$ for some $\lambda\ge0$.

\subsection{Ellipsoid}

Ellipsoidal uncertainty sets were first proposed in \cite{ben1998robust,ben1999robust} and stem from the observation that the iso-density locus of a multivariate normal distribution is an ellipse. We use an ellipsoid of the form $\cU^E = \{ \pmb{c} : (\pmb{c} - \pmb{\mu})^t\pmb{\Sigma}^{-1}(\pmb{c} - \pmb{\mu}) \le \lambda \}$ with size parameter $\lambda \ge 0$ that is centered on $\hat{\pmb{c}}$.
We create it using a normal distribution found as a maximum-likelihood fit. Recall that the best fit of a multivariate normal distribution $\mathcal{N}(\pmb{\mu},\pmb{\Sigma})$ with respect to data points $\pmb{c}^1,\ldots,\pmb{c}^N$ is given by
\[ \pmb{\mu} = \hat{\pmb{c}} = \frac{1}{N}(\pmb{c}^1 + \ldots + \pmb{c}^N)\]
and
\[ \pmb{\Sigma} = \frac{1}{N} \sum_{i\in[N]} (\pmb{c}^i - \pmb{\mu})(\pmb{c}^i - \pmb{\mu})^t \]
The resulting problem can then be formulated as 
%
\begin{align*}
\min\ & \hat{\pmb{c}}\pmb{x} + z \\
\text{s.t. } & z^2 \ge \lambda \left( \pmb{x}^t \pmb{\Sigma} \pmb{x}\right) \\ 
& \pmb{x}\in\X
\end{align*}
which is an integer second-order cone program (ISOCP), see \cite{ben1998robust} for details. Due to the convexity of the constraints, the problem can still be solved with little computational effort by standard solvers.



\subsection{Budgeted Uncertainty}

This approach was introduced in \cite{bertsimas2003robust}, and is based on 
interval uncertainty $\cU = \prod_{i\in[n]} \left[ \hat{c}_i, \overline{c}_i \right]$. To reduce the conservatism of this approach one assumes that only at most $\Gamma\in\{0,\ldots,n\}$ many values can be simultaneously higher than the midpoint $\hat{\pmb{c}}$. Formally,
\[\cU^{B} = \{ \pmb{c} : c_i = \hat{c}_i + (\overline{c}_i - \hat{c}_i)\delta_i \text{ for all } i\in[n],\ \pmb{0} \le \pmb{\delta}\le \pmb{1},\ \sum_{i\in[n]} \delta_i \le \Gamma\} \]
Using the dual of the inner worst-case problem, the following compact mixed-integer program can be found:
\begin{align*}
\min\ & \hat{\pmb{c}}\pmb{x} + \Gamma\pi + \sum_{i\in[n]} \rho_i \\
\text{s.t. } & \pi + \rho_i \ge (\overline{c}_i - \hat{c}_i)x_i & \forall i\in[n] \\
& \pi, \pmb{\rho} \ge 0 \\
& \pmb{x}\in\X
\end{align*}
This approach has the advantage that probability bounds can be found that compare favorably with those for ellipsoidal uncertainty \cite{bertsimas2004price}, while this problem also remains polynomially solvable by enumerating possible values for the $\pi$ variable. This means that $\mathcal{O}(n)$ many problems of the original type need to be solved. For these reasons, the budgeted uncertainty approach has been very popular in the literature.


\subsection{Permutohull}

The final two uncertainty sets we consider were proposed in \cite{bertsimas2009constructing}. The original inspiration comes from risk measures; the authors show that any so-called distortion risk measure leads to a polyhedral uncertainty set. A risk measure $\mu$ is a distortion risk measure if and only if there exists $\pmb{q}\in \{ \pmb{q}'\in\Delta^N : q_1 \ge \ldots q_N\}$, where $\Delta^N$ denotes the $N$-dimensional simplex such that
\[ \mu(\pmb{x}) = - \sum_{i\in[N]} q_i (\pmb{c}^{(i)} \pmb{x}) \]
where the sorting $(i)$ is chosen such that $\pmb{c}^{(1)} \pmb{x} \ge \ldots \ge \pmb{c}^{(N)} \pmb{x}$.

The conditional value at risk $CVaR_\alpha$ with $\alpha\in(0,1]$ is a well-known distortion risk measure. Intuitively, it is the expected value amongst the $\alpha$ worst outcomes. Using the matrix
\[ \pmb{Q}_N := \begin{pmatrix}
1 & \dots & \frac{1}{N-2} & \frac{1}{N-1} & \frac{1}{N} \\
\vdots & \vdots & \vdots & \vdots & \vdots \\
0 & 0 & \frac{1}{N-2} & \frac{1}{N-1} & \frac{1}{N} \\ 
0 & \dots & 0 & \frac{1}{N-1} & \frac{1}{N} \\
0 & \dots & 0 & 0 & \frac{1}{N}
\end{pmatrix}\]
the $j$th column of $Q_N$ induces the risk measure $CVaR_{j/N}$. The corresponding polyhedra are called the $\pmb{q}$-permutohull and defined as
\[ \Pi_{\pmb{q}}(\pmb{c}^1,\ldots,\pmb{c}^N) := conv\left( \left\{ \sum_{i\in[N]} q_{\sigma(i)} \pmb{c}^i : \sigma \in S_N\right\}\right) \]

To find the resulting robust problem, we first consider the worst-case problem for fixed $\pmb{x}\in\X$.
\begin{align*}
\max\ & \sum_{i,j\in[N]} q_i (\pmb{c}^j\pmb{x}) p_{ij} \\
\text{s.t. } & \sum_{i\in[N]} p_{ij} = 1 & \forall j\in[N] \\
& \sum_{j\in[N]} p_{ij} = 1 & \forall i\in[N] \\
& p_{ij} \ge 0 & \forall i,j\in[N]
\end{align*}
Dualising this problem then gives the robust counterpart
\begin{align*}
\min\ & \sum_{i\in[n]} (v_i + w_i) \\
\text{s.t. } & v_i + w_j \ge q_i (\pmb{c}^j\pmb{x}) & \forall i,j\in[N] \\
& \pmb{v},\pmb{w} \gtrless 0 \\
& \pmb{x}\in\X
\end{align*}
which is a mixed-integer program (note that this approach is actually the same as the ordered weighted averaging method, see \cite{chassein2015alternative}). The problem is NP-hard, as it contains the convex hull of $\{\pmb{c}^1,\ldots,\pmb{c}^N\}$ as a special case. 
%
Through the choice of $\pmb{q}$, there are $N$ possible sizes of this uncertainty.

\subsection{Symmetric Permutohull}

In the same setting as before, the symmetric permutohull was also introduced in \cite{bertsimas2009constructing}. By using the $\lfloor N/2 \rfloor +1$ columns of the matrix 
\[ \tilde{\pmb{Q}} := \frac{1}{N}\begin{pmatrix}
1 & 2 & 2 & \dots & 2 \\
1 & 1 & 2 & \dots & 2 \\
1 & 1 & 1 & \dots & 2 \\
\vdots & \vdots & \vdots  & \vdots  & \vdots  \\
1 & 1 & 1 & \dots & 0 \\
1 & 1 & 0 & \dots & 0 \\
1 & 0 & 0 & \dots & 0 \\
\end{pmatrix} \]
it was shown that the resulting polyhedra are symmetric with respect to $\hat{\pmb{c}}$. Note that these problems are also NP-hard, as $\tilde{\pmb{Q}}$ contains the min-max approach for $N=2$ as a special case.





\subsection{Summary of Uncertainty Sets}

\begin{figure}[htbp]
\centering
\begin{subfigure}[c]{0.45\textwidth}
        \includegraphics[width=\textwidth]{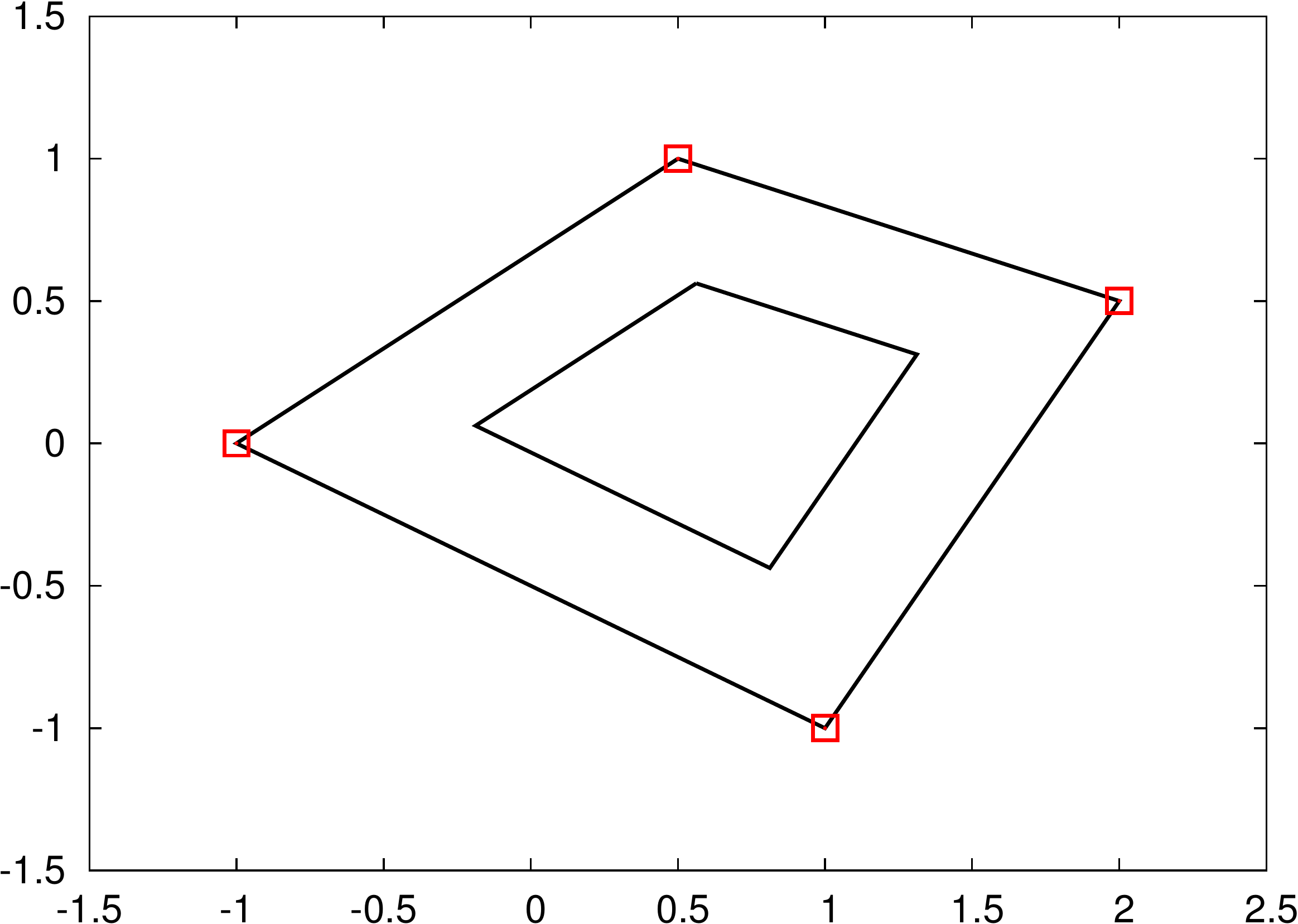}
        \caption{Convex hull with $\lambda=1$ and $\lambda=0.5$.}
        \label{graph1}
    \end{subfigure}
    \hfill
\begin{subfigure}[c]{0.45\textwidth}
        \includegraphics[width=\textwidth]{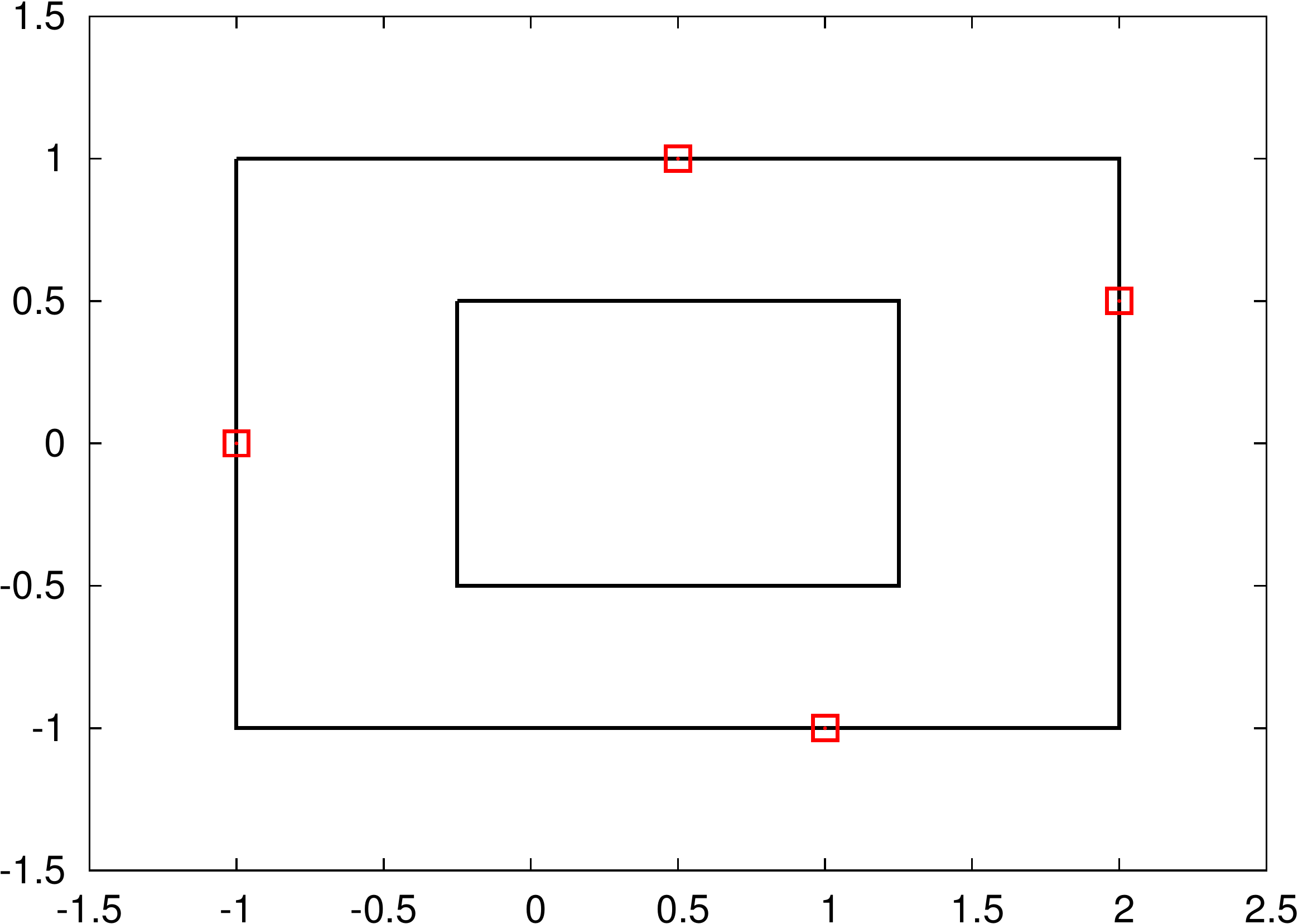}
        \caption{Intervals with $\lambda=1$ and $\lambda=0.5$.}
        \label{graph2}
    \end{subfigure}
\begin{subfigure}[c]{0.45\textwidth}
        \includegraphics[width=\textwidth]{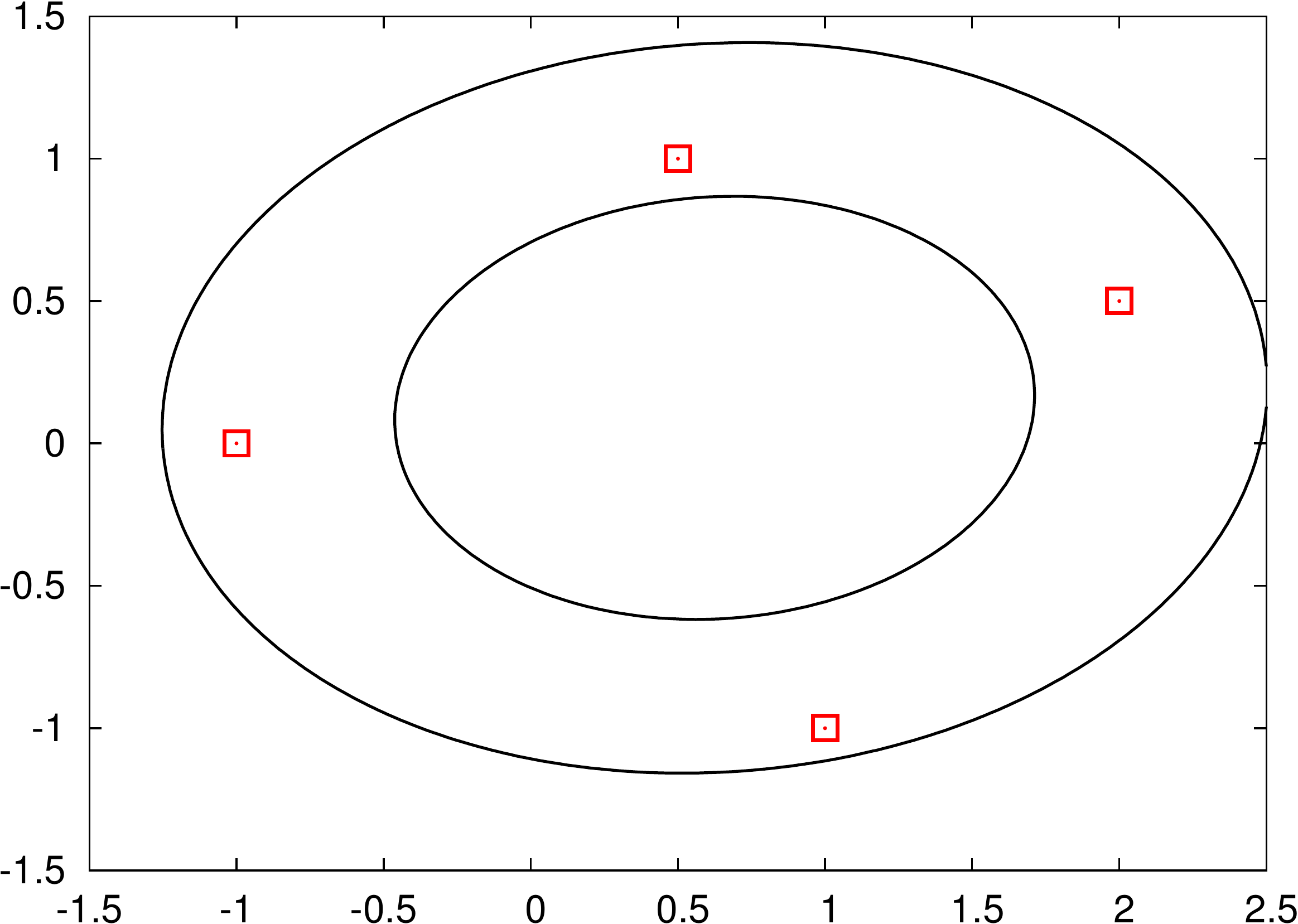}
        \caption{Ellipsoid with $\lambda=3$ and $\lambda=1$.}
        \label{graph3}
    \end{subfigure}   
    \hfill
\begin{subfigure}[c]{0.45\textwidth}
        \includegraphics[width=\textwidth]{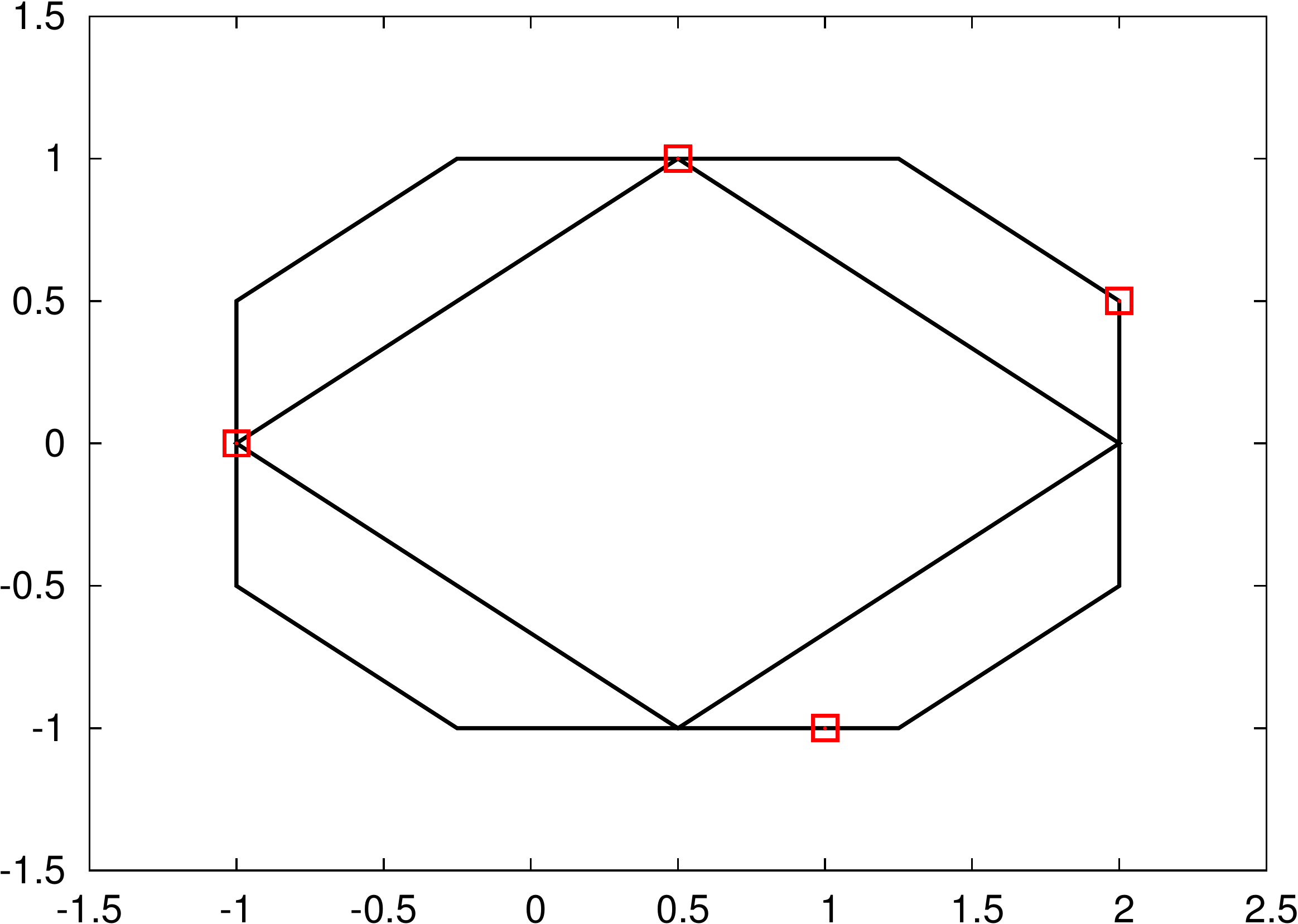}
        \caption{Budgeted uncertainty with $\lambda=1$, $\Gamma=1$ and $\lambda=1$, $\Gamma=1.5$.}
        \label{graph4}
    \end{subfigure}
    \begin{subfigure}[c]{0.45\textwidth}
        \includegraphics[width=\textwidth]{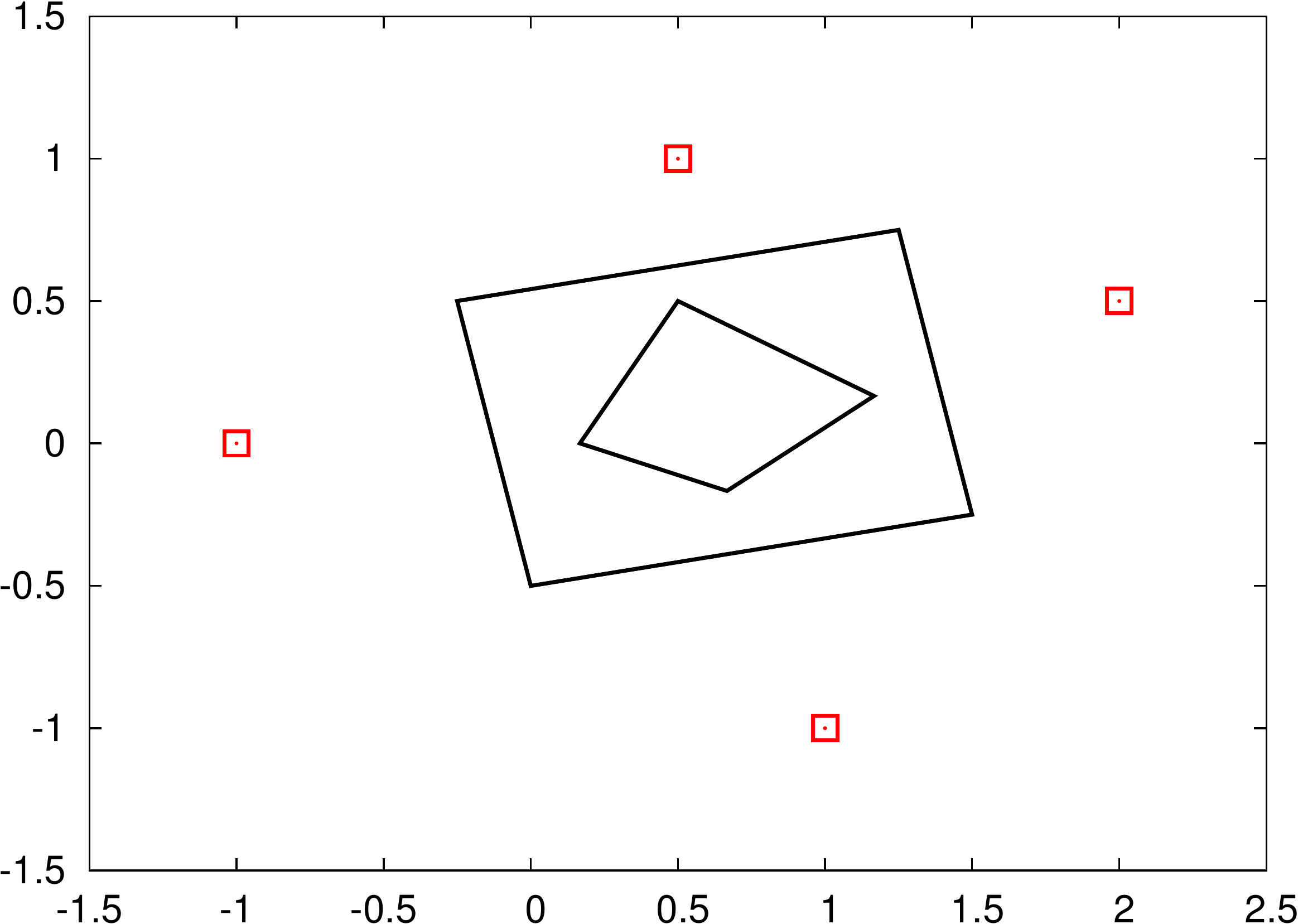}
        \caption{Permutohull uncertainty for $CVaR_{2/N}$ and $CVaR_{3/N}$ (i.e., $\pmb{q} = (\frac{1}{2},\frac{1}{2},0,0)$ and $\pmb{q}=(\frac{1}{3},\frac{1}{3},\frac{1}{3},0)$.}
        \label{graph5}
    \end{subfigure}     
    \hfill
    \begin{subfigure}[c]{0.45\textwidth}
        \includegraphics[width=\textwidth]{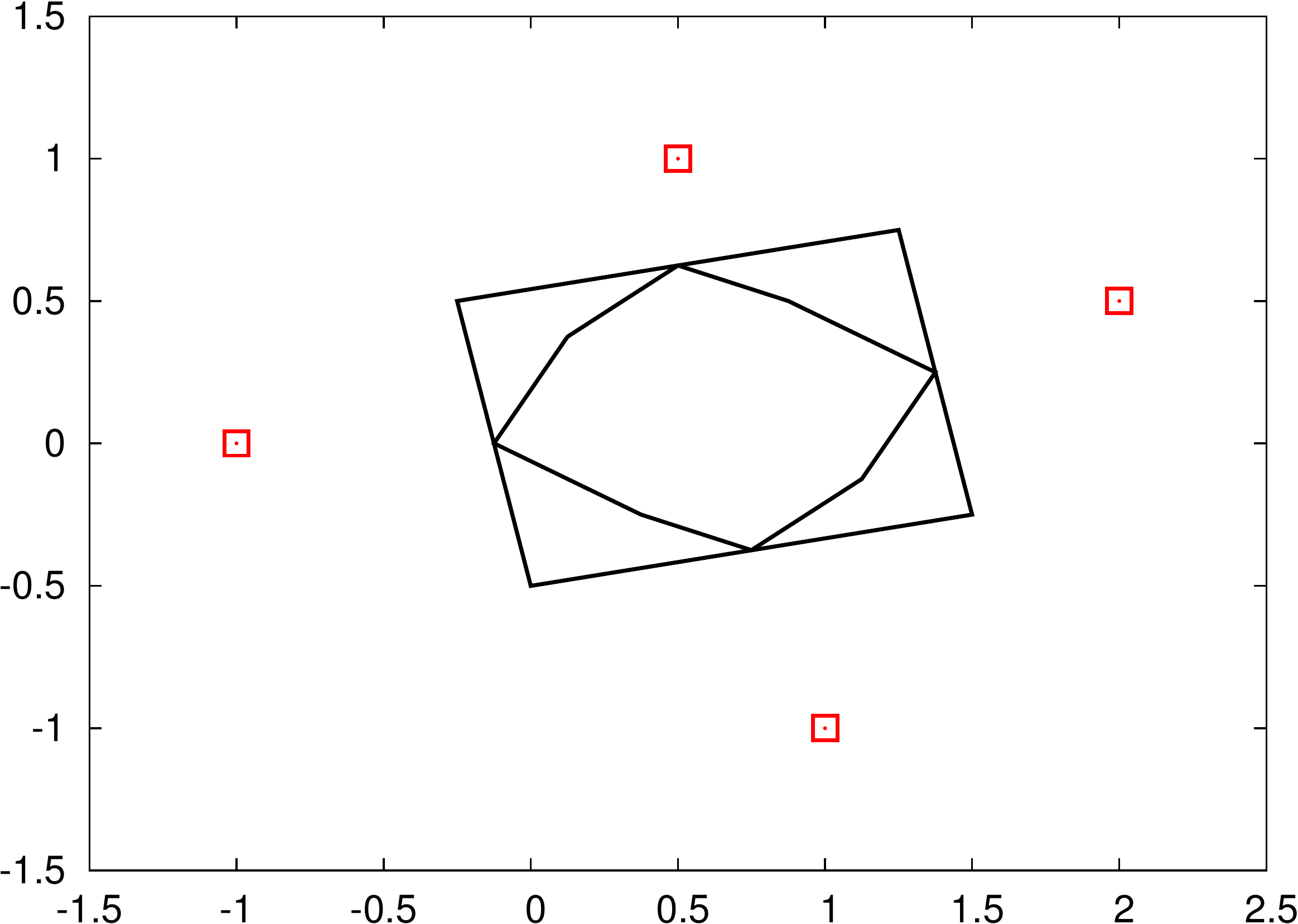}
        \caption{Symmetric permutohull uncertainty for $\pmb{q}=(\frac{1}{2},\frac{1}{2},0,0)$ and $\pmb{q}=(\frac{1}{2},\frac{1}{4},\frac{1}{4},0)$.}
        \label{graph6}
    \end{subfigure}  

\caption{Example uncertainty sets.}\label{fig1}
\end{figure}

In total we described six methods to generate uncertainty set $\cU$ based on the raw data $\cR$. Figure~\ref{fig1} illustrates these sets using a raw dataset with four observations (shown as red points). The complexity to solve the resulting robust models, as well as the type of program with the numbers of additional variables and constraints compared to the classic shortest path problem are shown in Table~\ref{tab1}.
\begin{table}[htb]
\centering
\begin{tabular}{r|c|c|c|c|c|c}
 & $\cU^{CH}$ & $\cU^{I}$ & $\cU^E$ & $\cU^{B}$ & $\cU^{PH}$ & $\cU^{SPH}$ \\
 \hline
Complexity & NPH & P & NPH & P & NPH & NPH \\
Model & IP & LP & ISOCP & MIP & MIP & MIP \\
Add. Const. & $N$ & 0 & 1 & $n+1$ & $N^2$ & $N^2$ \\
Add. Var. & 1 & 0 & 1 & $n$ & $2n$ & $2n$ \\
\end{tabular}
\caption{Uncertainty sets in this study.}\label{tab1}
\end{table}
While the robust model with budgeted uncertainty sets can be solved in polynomial time using combinatorial algorithms, we still used the MIP formulation for our experiments, as it was sufficiently fast.

\section{Real-World Experiments}
\label{sec:exp}
 
\subsection{Data Collection and Cleaning}

We used data provided by the City of Chicago\footnote{https://data.cityofchicago.org}, which provides a live traffic data interface. We recorded traffic updates in a 15-minute interval over a time horizon of 24 hours spanning Monday March 27th 2017 morning to Tuesday March 28th 2017 morning. A total of 98 data observations were thus used.

\begin{figure}[htbp]
\centering
\begin{subfigure}[c]{\textwidth}
        \includegraphics[width=\textwidth]{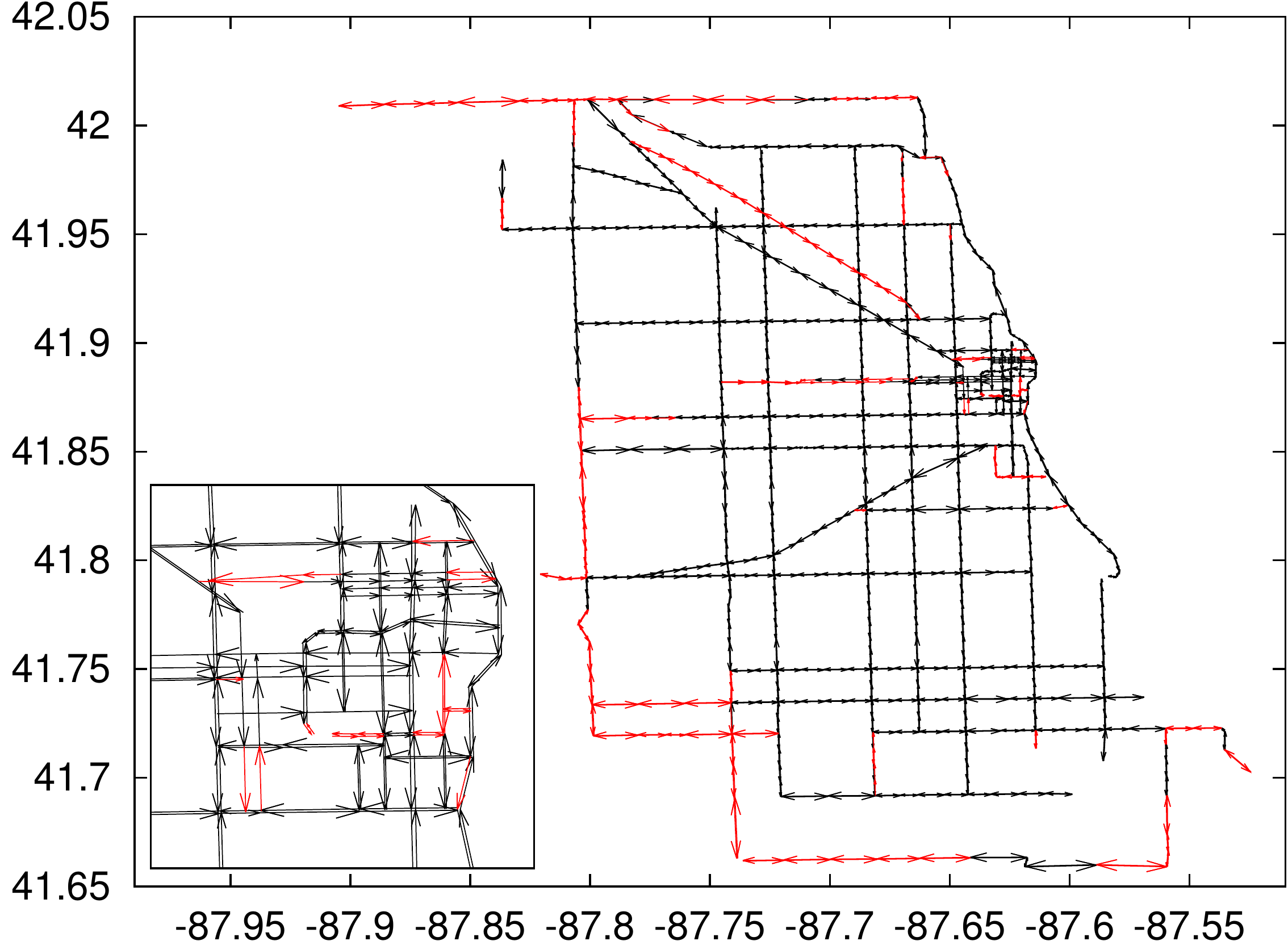}
        \caption{Raw segments with zoom-in for the city center. In red are segments without data.}
        \label{chicago1}
    \end{subfigure}
    \hfill
\begin{subfigure}[c]{\textwidth}
        \includegraphics[width=\textwidth]{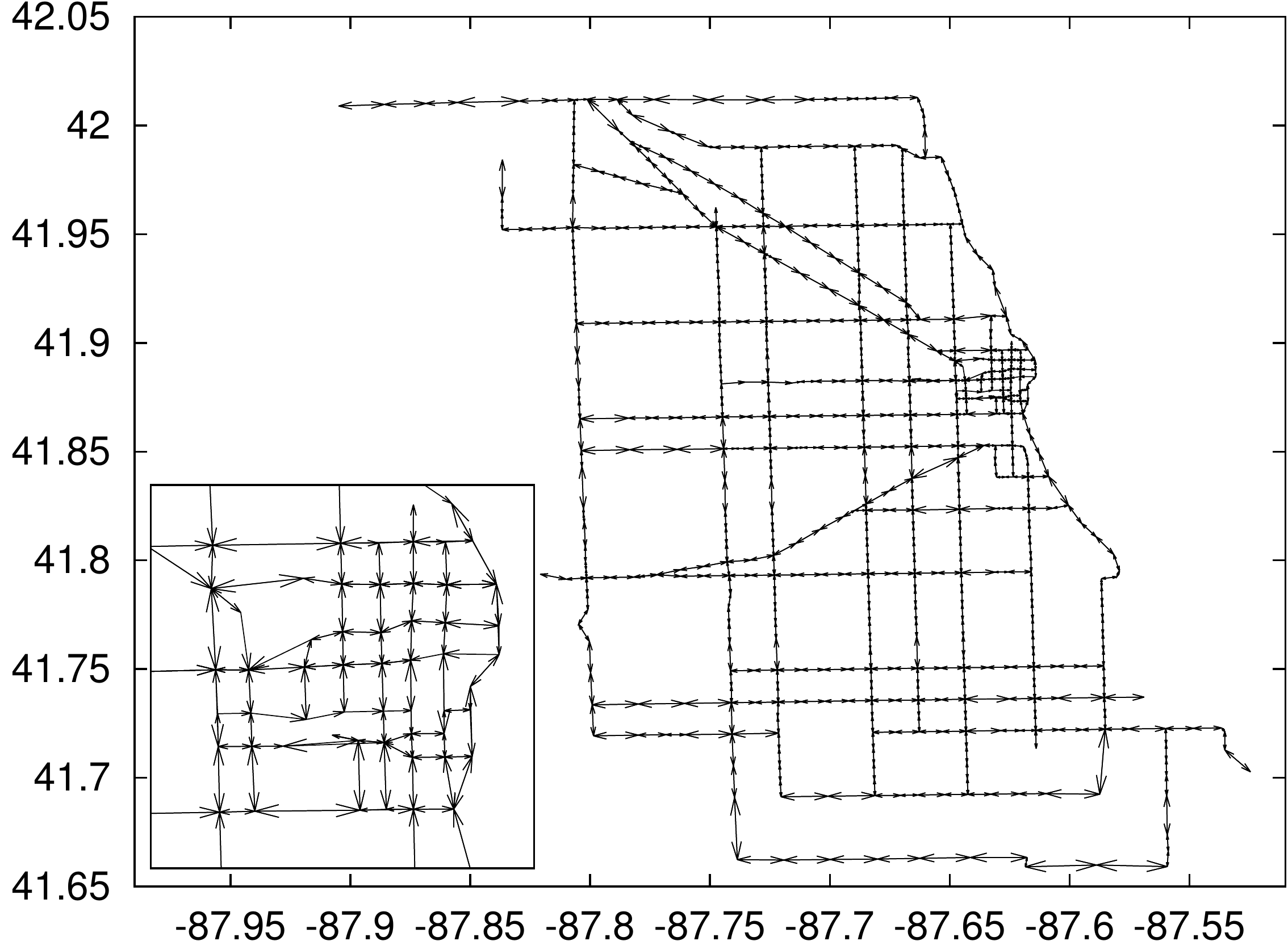}
        \caption{Resulting graph model with zoom-in for the city center.}
        \label{chicago2}
    \end{subfigure}
\caption{Chicago instance}
\end{figure}

Every observation contains the traffic speed for a subset of a total of 1,257 segments. For each segment the geographical position is available, see the resulting plot in Figure~\ref{chicago1} with a zoom-in for the city center. The complete travel speed data set contains a total of 54,295 observations. There were 1,027 segments where the data was recorded at least once of the 96 time points. Nearly for 88\% of the segments, speeds were recorded for at least 50 records with only 1\% (10 segments) where only one observation was recorded. We used linear interpolation to fill the missing records keeping in mind that data was collected over time. The data after removing missing records and filling missing values can be found at \url{www.lancaster.ac.uk/~goerigk/robust-sp-data.zip}.

As segments are purely geographical objects without structure, we needed to create a graph for our experiments. To this end, segments were split when they crossed or nearly crossed, and start- and end-points that were sufficiently close to each other were identified as the same node. The resulting graph is shown in Figure~\ref{chicago2}; note that this process slightly simplified the network, but kept its structure intact. The final graph contains 538 nodes and 1308 arcs.

\subsection{Setup}

Each uncertainty set is equipped with a size parameter. For each parameter we generated 20 possible values:
\begin{itemize}
\item For $\cU^{CH}$ and $\cU^I$, $\lambda \in \{ 0.1, 0.2, \ldots, 2\}$.
\item For $\cU^E$, $\lambda \in \{ 0.2, 0.4, \ldots, 4\}$.
\item For $\cU^B$, $\Gamma \in \{ 5, 10, \ldots, 100\}$.
\item For $\cU^{PH}$, we used columns $\pmb{q}_1, \pmb{q}_3, \ldots, \pmb{q}_{39}$.
\item For $\cU^{SPH}$, we used columns $\pmb{q}_1, \pmb{q}_2, \ldots, \pmb{q}_{20}$.
\end{itemize}
Each uncertainty set is generated using only every second scenario (i.e., 48 out of 96), but all 96 scenarios are then used to evaluate the solutions. Furthermore, we generated 200 random $s-t$ pairs uniformly, and used each of the $6\cdot 20$ methods on the same 200 pairs. Each of our 120 methods hence generates $200\cdot 96 = 19,200$ objective values.

It is highly non-trivial to assess the quality of these solutions, see \cite{chassein2016performance}. If one just uses the average objective value, as an example, then one could as well calculate the solution optimizing the average scenario case to find the best performance with respect to this measure. To find a balanced evaluation of all methods, we used four performance criteria:
\begin{itemize}
\item the average objective value,
\item the average of the worst-case objective value for each $s-t$ pair, and
\item the average value of the worst 5\% of objective values for each $s-t$ pair (as in the CVaR measure)
\end{itemize}
We also considered the average rank. To this end, we rank all 120 methods for each specific combination of $s-t$ pair and scenario. The best performing methods are ranked at 1, the second-best at 2 etc. We then take the average rank over all $19,200$ observations. However, this measure was strongly correlated with the average objective value and is therefore not presented.

For all experiments we used a computer with a 16-core Intel Xeon E5-2670 processor, running at 2.60 GHz with 20MB cache, and Ubuntu 12.04. Processes were pinned to one core. We used CPLEX v.12.6 to solve all problem formulations.

\subsection{Results}

\begin{figure}[htbp]
\centering
\begin{subfigure}[c]{\textwidth}
        \includegraphics[width=\textwidth]{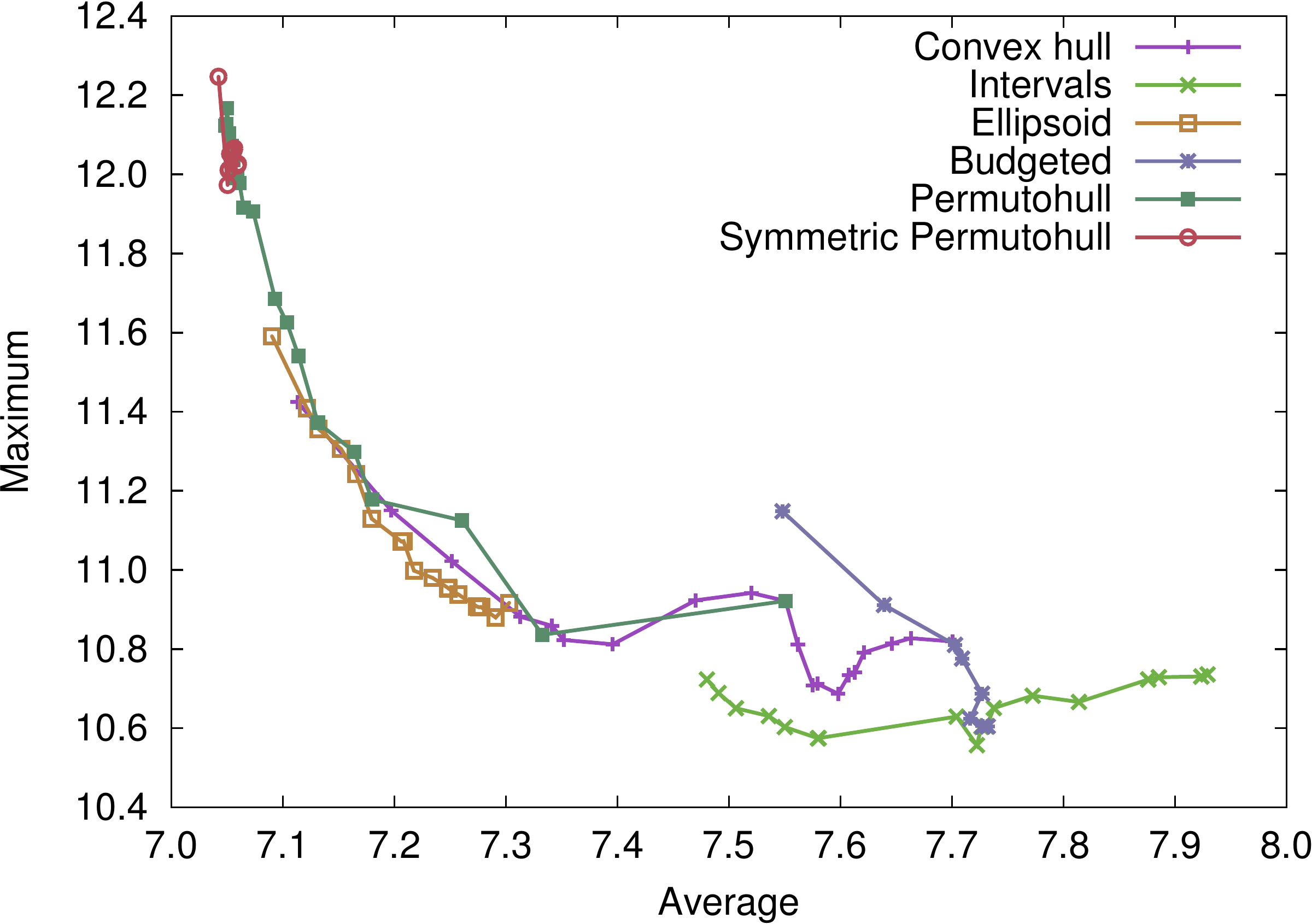}
        \caption{Average vs worst-case performance.}
        \label{exp1}
    \end{subfigure}
    \hfill
\begin{subfigure}[c]{\textwidth}
        \includegraphics[width=\textwidth]{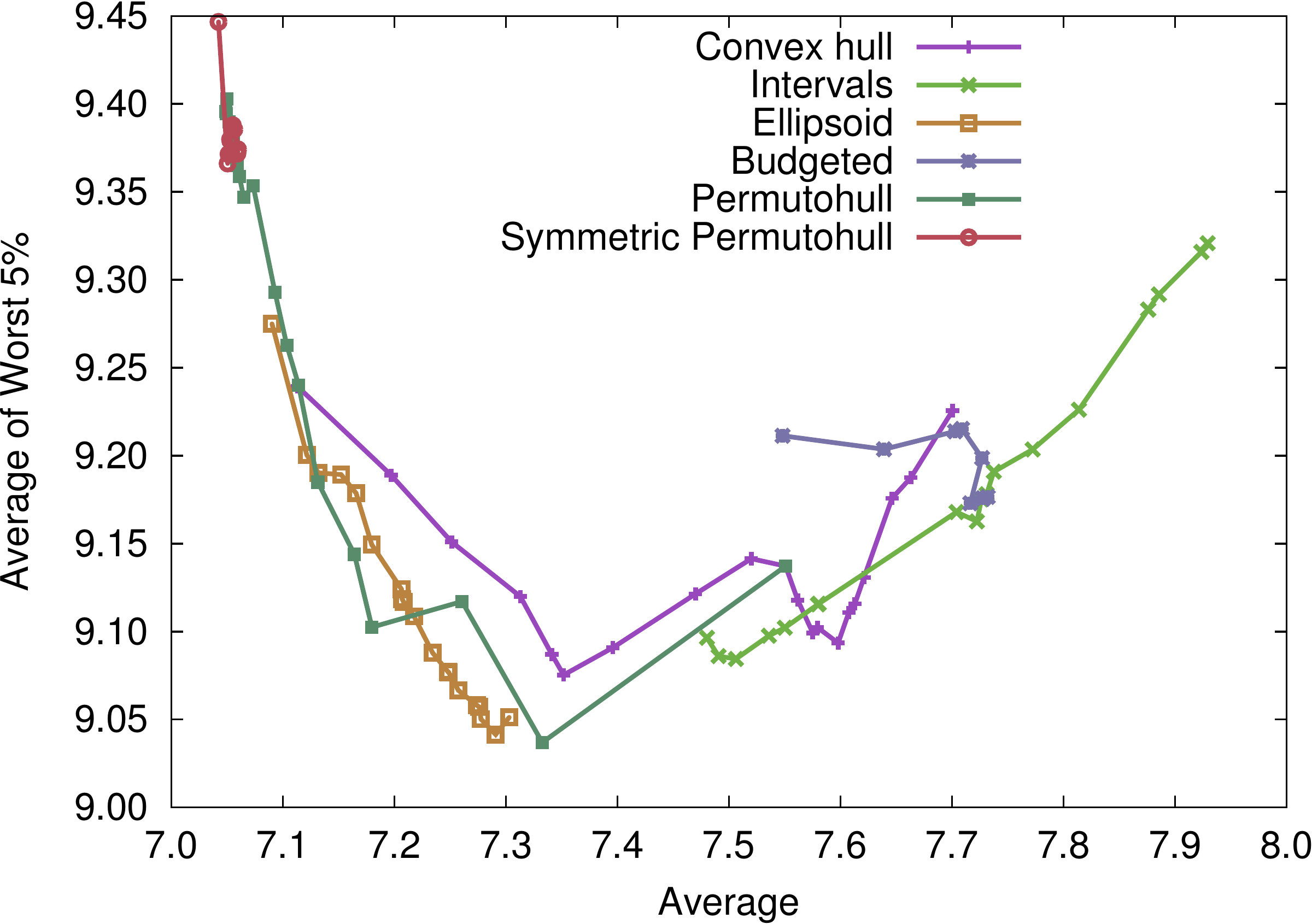}
        \caption{Average vs CVaR performance.}
        \label{exp2}
    \end{subfigure}
\caption{Performance results.}\label{results}
\end{figure}

\begin{figure}[htbp]
\centering
\includegraphics[width=\textwidth]{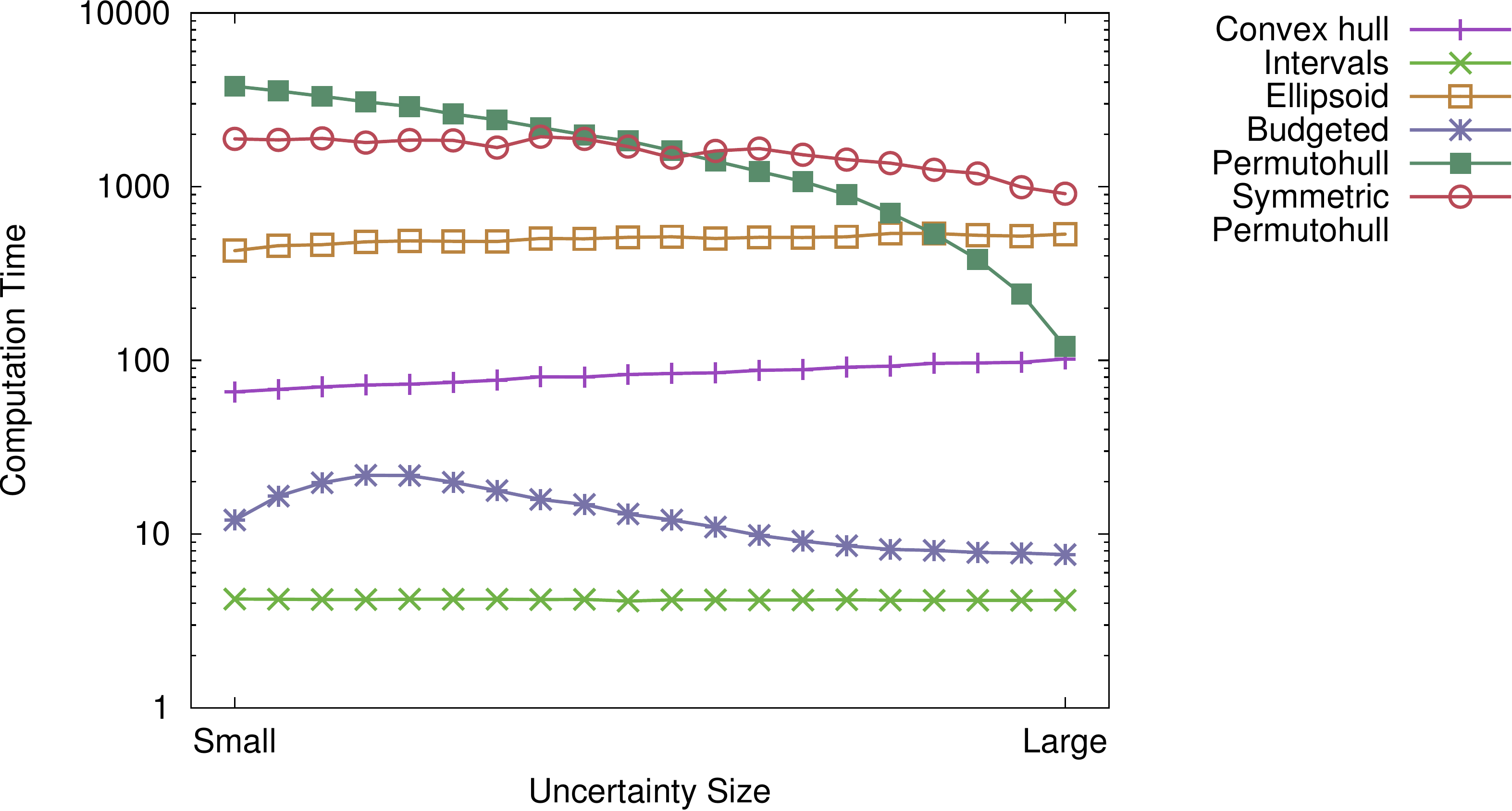}
\caption{Computation times.}\label{time}
\end{figure}

We present our findings in the two plots of Figure~\ref{results}. In each plot, the 20 parameter settings that belong to the same uncertainty set are connected by a line. They are complemented with Figure~\ref{time} showing the total computation times for the methods over all 200 shortest path calculations. 

The first plot in Figure~\ref{exp1} shows the trade-off between the average and the maximum objective value; the second plot in Figure~\ref{exp2} shows the trade-off between the average and the average of the 5\% worst objective values. Note that for all performance measures, smaller values indicate a better performance -- hence, good trade-off solutions should move from the top left to the bottom right of the plots. In general, the points corresponding to the parameter settings that give weight to the average performance can be on the left sides of the curves, while the more robust parameter settings are on the right sides, as would be expected.

We first discus Figure~\ref{exp1}. In general, we find that most concepts to indeed present a trade-off between average performance and robustness through their scaling parameter. The symmetric permutohull solutions have the best average performance, while interval solutions are the most robust. Interestingly, that even holds for interval solutions where the scaling parameter is very small. The budgeted uncertainty does not give a good trade-off between worst-case and average-case performance, which confirms previous results on artificial data \cite{chassein2016bicriteria}. Scaling interval uncertainty sets achieves better results than using budgeted uncertainty. Solutions generated with ellipsoidal uncertainty sets slightly outperform (dominate in the Pareto sense) solutions generated with permutohull uncertainty. We also note that methods that are computationally more expensive tend to achieve better average performance at the cost of decreases robustness. The simplest and cheapest method, interval uncertainty, gives the most robust solutions. Solutions using the convex hull of raw data tend to be outperformed by the approaches that process data.

We now consider the results presented in Figure~\ref{exp2}. Here the average is plotted against the average performance of the 5\% worst performing scenarios, averaged over all $s-t$ pairs. We note that for interval uncertainty, these two criteria are connected, with the best solutions for small parameter size dominating all solutions for larger parameter size. For the permutohull and the ellipsoidal uncertainty solutions, the order slightly changed with the former often dominating the latter. Permutohull solutions are designed to be efficient for the CVaR criterion, and the best-performing solution with respect to this aspect is indeed generated by this approach. However, also solutions with ellipsoidal, interval and convex hull uncertainty perform well.

Regarding computation times (see Figure~\ref{time}), note that the two polynomially solvable approaches are also the fastest when using Cplex; these computation times can be further improved using specialized algorithms. Using the convex hull is faster than using ellipsoids, which are in turn faster than using the symmetric permutohull. For the standard permutohull, the computation times are sensitive to the uncertainty size; if the $\pmb{q}$ vector that is used in the model has only few entries, then computation times are smaller. This is in line with the intuition that the problem becomes easier if less scenarios need to be considered.

To summarize our findings in our experiment on the robust shortest path problem with real-world data:
\begin{itemize}
\item Convex hull solutions are amongst the more robust solutions, but tend to be outperformed by the other approaches.
\item Interval solutions perform bad on average, but are the most robust. Especially when the scaling is small they can give a decent trade-off, and are easy and fast to compute.
\item Ellipsoidal uncertainty solutions have very good overall performance and represent a large part of the non-dominated points in our results.
\item We do not encourage the use of budgeted uncertainty for robust shortest path problems. Scaling interval uncertainty sets gives better results and is easier to use and to solve.
\item Permutohull solutions offer good trade-off solutions, whereas symmetric permutohull solutions tend to be less robust, but provide an excellent average performance. These methods are also require most computational effort to find.
\end{itemize}
In the light of these findings, the interval and discrete (=convex hull) uncertainty sets that are widely used in robust combinatorial optimization do warrant research attention, as they may not produce the best solutions, but are relatively fast to solve. However, permutohull and ellipsoidal uncertainty tend to produce solutions with a better trade-off, while being computationally more challenging. The algorithmic research for robust shortest path problems with such structure should therefore become a future focus.

\section{Conclusion}
\label{sec:conclusion} 
 
In this paper wo constructed uncertainty sets for the robust shortest path problem using real-world traffic observations for the City of Chicago. We evaluated the model suitability of these sets by finding the resulting robust paths, and comparing their performance using different performance indicators.

Naturally, conclusions can only be drawn within the reach of the available data. In our setting we considered solutions that are robust with respect to all possible travel times within a day. A use-case would be that a path needs to be computed for a specific day, but the precise hour is not known. Using different sets of observations will result in solutions that are different in another sense, e.g., one could use observations over different days during the morning rush hours, or observations that span work days and a weekend. It is possible that these sets will provide different structure.
 
Finally, we have observed that using ellipsoidal uncertainty sets provides high-quality solutions with less computational effort than for the permutohull. If one uses only the diagonal entries of the matrix $\pmb{\Sigma}$, then one ignores the data correlation in the network. For the resulting problem specialized algorithms exist, see, e.g. \cite{nikolova2009high}. In additional experiments we found that even by using Cplex, computation times were considerably reduced when only using the diagonal entries of $\pmb{\Sigma}$, but the solution quality remained roughly the same.

\newcommand{\etalchar}[1]{$^{#1}$}

\end{document}